\documentclass[11pt]{amsart}
\usepackage{epsfig}
\usepackage[all]{xy}
\begin{document}
\title{Vari\'et\'es projectives complexes dont l'\'eclat\'ee en un point est
de Fano}
\author{Laurent BONAVERO, Fr\'ed\'eric CAMPANA et Jaros{\l}aw A. 
WI\'SNIEWSKI}
\date{Juillet 2001, version r\'evis\'ee de math.AG/0106047}
\maketitle
\noindent
\def\restriction{\string |}
\newcommand{\pp}{\rm ppcm}
\newcommand{\pg}{\rm pgcd}
\newcommand{\Ker}{\rm Ker}
\newcommand{\C}{{\mathbb C}}
\newcommand{\Q}{{\mathbb Q}}
\newcommand{\GL}{\rm GL}
\newcommand{\SL}{\rm SL}
\newcommand{\diag}{\rm diag}

\def\refname{R\'ef\'erences}
\def\finpreuve
{\hskip 3pt \vrule height6pt width6pt depth 0pt}

\newtheorem{theo}{Th\'eor\`eme}
\newtheorem{prop}{Proposition}
\newtheorem{lemm}{Lemme}
\newtheorem{lemmf}{Lemme fondamental}
\newtheorem{defi}{D\'efinition}
\newtheorem{exo}{Exercice}
\newtheorem{rem}{Remarque}
\newtheorem{cor}{Corollaire}
\newcommand{\CC}{{\mathbb C}}
\newcommand{\ZZ}{{\mathbb Z}}
\newcommand{\RR}{{\mathbb R}}
\newcommand{\QQ}{{\mathbb Q}}
\newcommand{\FF}{{\mathbb F}}
\newcommand{\PP}{{\mathbb P}}
\newcommand{\codim}{\operatorname{codim}}
\newcommand{\Ho}{\operatorname{Hom}}
\newcommand{\Pic}{\operatorname{Pic}}
\newcommand{\NE}{\operatorname{NE}}
\newcommand{\Nun}{\operatorname{N}}
\newcommand{\card}{\operatorname{card}}
\newcommand{\Hilb}{\operatorname{Hilb}}
\newcommand{\mult}{\operatorname{mult}}
\newcommand{\vol}{\operatorname{vol}}
\newcommand{\divi}{\operatorname{div}}
\newcommand{\pr}{\operatorname{pr}}
\newcommand{\con}{\operatorname{cont}}
\newcommand{\Spec}{\operatorname{Spec}}
\newcommand{\ima}{\operatorname{Im}}

\newcounter{subsub}[subsection]
\def\thesubsub{\thesubsection .\arabic{subsub}}
\def\subsub#1{\addtocounter{subsub}{1}\par\vspace{3mm}
\noindent{\bf \thesubsub ~ #1 }\par\vspace{2mm}}
\def\coker{\mathop{\rm coker}\nolimits}
\def\pr{\mathop{\rm pr}\nolimits}
\def\im{\mathop{\rm Im}\nolimits}
\def\hfl#1#2{\smash{\mathop{\hbox to 12mm{\rightarrowfill}}
\limits^{\scriptstyle#1}_{\scriptstyle#2}}}
\def\vfl#1#2{\llap{$\scriptstyle #1$}\big\downarrow
\big\uparrow
\rlap{$\scriptstyle #2$}}
\def\diagram#1{\def\normalbaselines{\baselineskip=0pt
\lineskip=10pt\lineskiplimit=1pt}   \matrix{#1}}
\def\limind{\mathop{\oalign{lim\cr
\hidewidth$\longrightarrow$\hidewidth\cr}}}

\long\def\InsertFig#1 #2 #3 #4\EndFig{
\hbox{\hskip #1 mm$\vbox to #2 mm{\vfil\includegraphics{#3}}#4$}}
\long\def\LabelTeX#1 #2 #3\ELTX{\rlap{\kern#1mm\raise#2mm\hbox{#3}}}

{\let\thefootnote\relax
\footnote{
\textrm{Mots cl\'es : vari\'et\'es de Fano, \'eclatements, 
th\'eorie de Mori}.
\textrm{Classification~A.M.S. : 14J45, 14E30.} 
}}

{\let\thefootnote\relax
\footnote{
\textrm{J.A.\ W. b\'en\'eficie du soutien financier du
contrat 2P03A02216
du KBN Polonais.}
}}

{\bf R\'esum\'e : } Nous classifions les vari\'et\'es projectives
complexes $X$ pour lesquelles il existe un point
$a$ tel que l'\'eclatement de $X$ en $a$ soit une vari\'et\'e de Fano. 

\bigskip

{\bf Abstract : } We classify complex projective manifolds $X$
for which there exists a point $a$ such that
the blow-up of $X$ at $a$ is Fano.

\bigskip

\section*{Introduction}
La g\'eom\'etrie birationnelle des vari\'et\'es de Fano,
{\em i.e.} des vari\'et\'es dont le fibr\'e anti-canonique est ample, 
intervient de fa\c con essentielle en th\'eorie de Mori
et le comportement des vari\'et\'es de Fano par \'eclatements
ou contractions lisses 
est une question naturelle (\'etudi\'ee en particulier dans 
\cite{Wis91}) \`a
la suite d'une part du th\'eor\`eme de factorisation
des applications birationnelles \cite{Aal99} et d'autre
part de l'\'etude torique r\'ecente
men\'ee par Sato \cite{Sat00} (voir aussi \cite{Cas01}).
Nous d\'emontrons dans cette Note le r\'esultat suivant (voir
aussi \cite{Kol99} exercise V.3.7.10)~:

\begin{theo}\label{main}
Soit $X$ une vari\'et\'e projective complexe lisse et connexe
de dimension $n \geq 3$. Soient $a \in X$
et $\pi_a~: \tilde X \to X$ l'\'eclatement de $X$ en $a$.
Alors $\tilde X$ est une vari\'et\'e de Fano si 
et seulement si~: 
\begin{enumerate} 
\item [(i)] soit $X$ est isomorphe \`a l'espace projectif $\PP^n$
et $a$ est quelconque dans $X$,
 \item [(ii)] soit $X$ est isomorphe \`a la quadrique ${\mathcal Q}_n$
et $a$ est quelconque dans $X$,
\item [(iii)] soit $X$ est isomorphe \`a la vari\'et\'e $V_d$ 
obtenue en \'eclatant 
$\PP ^n$ le long d'une sous-vari\'et\'e lisse $Y$ de dimension $n-2$,
de degr\'e $d$ avec
$1 \leq d \leq n$ 
contenue dans un hyperplan $H$ et $a \notin H$.
\end{enumerate}

\end{theo}

Remarquons que si l'hypoth\`ese du th\'eor\`eme \ref{main}
est satisfaite pour une vari\'et\'e $X$, alors $X$ est 
n\'ec\'essairement une vari\'et\'e de Fano.
Notons que le premier auteur a obtenu la version torique
de ce r\'esultat en classifiant plus g\'en\'eralement
les vari\'et\'es toriques de Fano
de dimension $n\geq 3$  
contenant un diviseur torique isomorphe \`a $\PP ^{n-1}$ \cite{Bon00}.
En dimension trois, le th\'eor\`eme \ref{main} peut se d\'emontrer
directement \`a l'aide de la classification des vari\'et\'es
de Fano de dimension trois \cite{MM81}.

Le th\'eor\`eme \ref{main} poss\`ede le corollaire imm\'ediat suivant
qui caract\'erise la quadrique en termes d'\'eclatements
ponctuels~:

\begin{cor}
Soit $X$ une vari\'et\'e projective complexe lisse et connexe
de dimension $n \geq 3$. Soient $a$ et $b$ deux points 
distincts de $X$
et $\pi_{a,b}~: \tilde X \to X$ l'\'eclatement de $X$ 
de centre $\{a,b\}$.
Alors $\tilde X$ est une vari\'et\'e de Fano si 
et seulement si $X$ est isomorphe \`a la quadrique ${\mathcal Q}_n
\subset \PP^{n+1}$
et $a$ et $b$ ne sont pas sur une m\^eme droite de ${\mathcal Q}_n$.
\end{cor}

\noindent {\bf D\'emonstration.} En effet, seule la vari\'et\'e
$V_2$ du th\'eor\`eme \ref{main} poss\`ede un    
diviseur isomorphe \`a $\PP ^{n-1}$ avec fibr\'e normal
${\mathcal O}_{\PP^{n-1}}(-1)$ 
(donn\'ee par la transform\'ee stricte
de l'hyperplan $H$) et la contraction sur un point de ce dernier
est ${\mathcal Q}_n$.\finpreuve

\section{Pr\'eliminaires}
Si $X$ est une vari\'et\'e projective, on note 
$$\Nun_1(X) = \{ \sum_i a_i C_i \, | \, a_i \in \QQ \, , \, 
C_i \, \, \mbox{courbe irr\'eductible de}\,  \, X \}/\equiv $$
o\`u $\equiv$ d\'esigne l'\'equivalence num\'erique. 
Le c\^one de Mori, ou c\^one des courbes effectives,
est le sous-c\^one de $\Nun_1(X)$ d\'efini par
$$ \NE (X) = \{ Z\in \Nun_1(X) \, | \, Z \equiv 
\sum_i a_i C_i \, , \, a_i\geq 0\}.$$
Si $X$ est non singuli\`ere de Fano,
alors $\NE (X)$ est poly\'edral, 
engendr\'e par des classes de
courbes rationnelles.
De plus, pour
toute ar\^ete $R$ de $\NE (X)$, il y a
une contraction $\varphi_R~: X \to Y$ de $X$ sur une vari\'et\'e 
projective \'eventuellement singuli\`ere $Y$ telle que 
les courbes irr\'eductibles contract\'ees par $\varphi_R$ sont exactement
celles dont la classe dans $\Nun_1(X)$ appartient \`a
$R$.

Terminons ces pr\'eliminaires en montrant que les vari\'et\'es
$V_d$ satisfont bien l'hypoth\`ese du th\'eor\`eme \ref{main}, 
les paragraphes suivants seront consacr\'es
\`a la r\'eciproque. 

\begin{lemm}\label{direct}
Soit $W_d$ la vari\'et\'e  
obtenue en \'eclatant 
$\PP ^n$ le long d'une sous-vari\'et\'e lisse $Y$ de dimension $n-2$,
de degr\'e $d$ contenue dans un hyperplan $H$ 
de $\PP^n$ et en un point $a$ n'appartenant pas \`a $H$.
Alors $W_d$ est de Fano
si et seulement si $1\leq d \leq n$.   
\end{lemm}

\noindent {\bf D\'emonstration.}
Reprenons les notations du th\'eor\`eme \ref{main} (iii)
et notons $\mu~: V_d \to \PP^n$ l'\'eclatement d\'efinissant $V_d$.
Si $H' \subset V_d$ est la transform\'ee stricte de $H$,
le c\^one $\NE (V_d)$ est le c\^one ferm\'e de $\Nun_1(V_d)$
engendr\'e par la classe d'une fibre non triviale $f$ de $\mu$
et d'une droite $l$ contenue dans $H'$.
Comme $-K_{V_d}\cdot f =1$ et $-K_{V_d}\cdot l = n+1-d$,
on en d\'eduit, d'apr\`es
le crit\`ere de Kleiman \cite{Kle66},
que $V_d$ est de Fano si et seulement si $1\leq d \leq n$,
condition que l'on suppose d\'esormais satisfaite.
De l\`a, soient $E \subset V_d$ le diviseur exceptionnel
de $\mu$ et
$\tilde E \subset W_d$ celui de $\pi_a$.
Si $C$ est une courbe de $W_d$ non contenue dans $\tilde E$,
\'ecrivons $[(\pi_a)_* (C)] = \delta [l] + \beta [f]$
dans $\NE (V_d)$ o\`u $\delta$ et $\beta$ sont des rationnels
positifs. Alors 
$-K_{W_d}\cdot C = (n+1)\delta - (n-1) \tilde E \cdot C
- E \cdot ((\pi_a)_* (C))$, et comme
$\tilde E \cdot C \leq \delta$ et  $E \cdot ((\pi_a)_* (C)) \leq \delta$
par le th\'eor\`eme de Bezout, on en d\'eduit
que $-K_{W_d}\cdot C \geq \delta$ 
et \`a nouveau par le crit\`ere de Kleiman
que $W_d$ est de Fano. \finpreuve

\section{R\'esultats pr\'eparatoires}

Dans ce travail, si $V$ est un espace vectoriel,
$\PP (V)$ d\'esigne l'espace projectif des {\em droites} de $V$.

\begin{lemm}\label{existcont}
Soient $\tilde X$ une vari\'et\'e de Fano 
de dimension $n$ et $\tilde D \simeq \PP ^{n-1}$
un diviseur de $\tilde X$.
Alors il existe une contraction extr\'emale
$\varphi ~: \tilde X \to X'$ d'ar\^ete $R = \QQ ^+ [\tilde F]$
avec $\tilde F \cdot \tilde D > 0$ o\`u
$\tilde F$ est une courbe rationnelle engendrant $R$.
\end{lemm} 

\noindent {\bf D\'emonstration.} Soit $\tilde C$ une courbe
telle que $\tilde C \cdot \tilde D > 0$.
Cette courbe est nu\-m\'e\-ri\-quement combinaison lin\'eaire \`a coefficients
rationnels strictements positifs de courbes extr\'emales. L'une d'entre
elles, not\'ee $\tilde F$ satisfait $\tilde F \cdot \tilde D > 0$.
\finpreuve  

\medskip

La proposition suivante pr\'ecise la structure 
de la contraction extr\'emale donn\'ee par le lemme~\ref{existcont}.

\begin{prop}\label{Step1} Soient $\tilde X$ une vari\'et\'e de Fano 
de dimension $n \geq 3$, $\tilde D \simeq \PP ^{n-1}$
un diviseur de $\tilde X$ de fibr\'e normal 
$N_{\tilde D / \tilde X} \simeq {\mathcal O}_{\PP ^{n-1}}(-1)$ 
et $\varphi ~: \tilde X \to X'$
une contraction extr\'emale
d'ar\^ete $R = \QQ ^+ [\tilde F]$
avec $\tilde F \cdot \tilde D > 0$. 
Alors $\varphi ~: \tilde X \to X'$ est
de l'un des deux types suivants~:
\begin{enumerate}
\item [a)] la fibration lisse en $\PP^1$ sur $\PP ^{n-1}$ donn\'ee par
$\varphi~:\tilde X \simeq \PP ({\mathcal O}_{\PP ^{n-1}} \oplus
{\mathcal O}_{\PP ^{n-1}}(-1)) \to X' \simeq \PP ^{n-1}$ 
et $\tilde D$ est le diviseur
$\PP ({\mathcal O}_{\PP ^{n-1}})$~;
\item [b)] une contraction lisse de centre $Y' \subset X'$ lisse 
de codimension deux ({\em i.e.} $X'$ est une vari\'et\'e projective lisse et 
$\tilde X$ est obtenue en \'eclatant $X'$ le long de $Y'$).
De plus, $\varphi _{|\tilde D}~:
\tilde D \to D':= \varphi (\tilde D)$ est un isomorphisme, 
$Y' \subset D'$ est un diviseur de degr\'e $d'+1$ 
o\`u $d'$ est le degr\'e du fibr\'e normal $N_{D'/X'}$ et
$X'$ est de Fano. 
\end{enumerate}

\end{prop}

\noindent {\bf D\'emonstration.}\footnote{Nous remercions
le referee anonyme qui nous a mentionn\'e les erreurs contenues dans
une premi\`ere version de ce travail.}
Montrons tout d'abord que les fibres non triviales de
$\varphi$ sont de dimension un. En effet, si
$f$ est une fibre non triviale de $\varphi$,
et si $\dim f \geq 2$, alors $\dim (f \cap \tilde D) \geq 1$
et l'ar\^ete $R$ contient donc une courbe contenue dans
$\tilde D$, ceci est absurde puisque $R \cdot \tilde D > 0$
alors que $N_{\tilde D / \tilde X} \simeq {\mathcal O}_{\PP ^{n-1}}(-1)$.
Il d\'ecoule d'un r\'esultat de Ando \cite{And85}
que $X'$ est lisse et que $\varphi~: \tilde X \to X'$
est de l'un des trois types suivants~: 
\begin{enumerate}
\item [(i)] une fibration lisse en $\PP ^1$,
\item [(ii)] une fibration en coniques (ce qui signifie
qu'il y a effectivement des fibres singuli\`eres),  
\item [(iii)] une contraction lisse de centre $Y' \subset X'$ lisse 
de codimension deux.
\end{enumerate}

Dans les cas (i) et (ii), $\varphi _{|\tilde D}~:
\tilde D \to X'$ est finie, et comme $\tilde D \simeq \PP^{n-1}$,
il s'ensuit que $X' \simeq \PP ^{n-1}$ d'apr\`es \cite{Laz83}
et par suite que $\rho(\tilde X)=2$.
\begin{enumerate}
\item [$\bullet$]
Dans le cas (i), soit $l$ une droite de $X'$
et $S_l$ la surface r\'egl\'ee $S_l = \varphi ^{-1}(l)$
($S_l$ est une surface de Hirzebruch ${\mathbf F}_a$). 
Comme $\tilde D$ est exceptionnel dans $\tilde X$, 
$\tilde D \cap S_l$ l'est dans $S_l$.
Il s'ensuit que $\tilde D$ est une section de $\varphi$
et de l\`a $\tilde X \simeq \PP ({\mathcal O}_{\PP ^{n-1}} \oplus
{\mathcal O}_{\PP ^{n-1}}(-1))$.

\item [$\bullet$] Excluons le cas (ii). Dans ce cas, 
il y a une courbe $C$ (composante d'une fibre singuli\`ere)
sur l'ar\^ete $R$ telle
que $-K_{\tilde X} \cdot C =1$, 
${\mathcal E} := \varphi_*(-K_{\tilde X})$
est un
fibr\'e vectoriel de rang $3$ sur $X'$ 
(voir par exemple \cite{ABW93} Theorem~B) et  
il y a un diagramme
commutatif~:

\centerline{
\xymatrix{ \tilde{X}  \ar@{^{(}->}[r]  \ar[rd]^{\varphi}&
   \PP({\mathcal E}^*) \ar[d]^{\mu} \\
& X'  
}
}
\noindent 
Cette situation est exclue par le lemme \ref{conique} ci-dessous
car $N_{\tilde D / \tilde X} \simeq {\mathcal O}_{\PP ^{n-1}}(-1)$.

\item [$\bullet$] Dans le cas (iii), notons $\tilde E$ le
diviseur exceptionnel de $\varphi$.
Comme $R \cdot \tilde D > 0$, chaque fibre non triviale 
de $\varphi$ rencontre $\tilde D$. Puisque
$\varphi _{|\tilde E}~: \tilde E \to Y'$ est 
une fibration lisse en $\PP ^1$ et que
$\tilde D \cap \tilde E$ est exceptionnel dans $\tilde E$,
il s'ensuit que $\tilde D \cap \tilde E$ est une section 
({\em a priori} non n\'ecessairement r\'eduite) de 
$\varphi _{|\tilde E}$.

{\em Admettons un instant} que cette section est r\'eduite
(autrement dit, $\tilde D$ et $\tilde E$ s'in\-ter\-sec\-tent
transversalement).
De l\`a, $\varphi _{|\tilde D}~: \tilde D \to D':=\varphi (\tilde D)$
est un isomorphisme. Donc
$Y' \subset D' \simeq \PP^{n-1}$ comme diviseur
de degr\'e $d'+1$ si $d'\geq 0$ d\'esigne le degr\'e
du fibr\'e normal $N_{D'/X'}$. Il s'ensuit aussi que $X'$ est de Fano. 

Il nous reste \`a v\'erifier que si  
$Y_0 := (\tilde D \cap \tilde E)_{red}$, la
multiplicit\'e d'intersection de
$\tilde E$ et $\tilde D$ le long de $Y_0$ est \'egale \`a un. 
Autrement dit, \'ecrivons 
$\tilde D _{|\tilde {E}}=mY_0$ et montrons que $m=1$.
Soit $N$ le fibr\'e normal de $Y_0$ dans $\tilde E$. 
Dans $\Pic (Y_0)$, on a $\tilde D _{|Y_0}=mN$.
Or, pour $n\geq 4$, l'application de restriction 
$\Pic (\tilde D) \simeq \ZZ \to \Pic (Y_0)$
est injective et son conoyau est sans torsion d'apr\`es le th\'eor\`eme de
Lefschetz (voir par exemple 
\cite{BS95} page 51). Comme $\tilde D _{|\tilde {D}} \simeq  
{\mathcal O}_{\PP ^{n-1}}(-1)$ est un g\'en\'erateur de $\Pic (\tilde D)$,
il s'ensuit que $m=1$ (pour $n=3$, la proposition \ref{Step1}
comme le th\'eor\`eme \ref{main} d\'ecoulent par exemple de \cite{MM81}).

\end{enumerate}
\finpreuve

\begin{lemm}\label{conique}
Soit $Z$ une vari\'et\'e lisse de dimension $n\geq 3$ 
telle que~:

\centerline{
\xymatrix{ Z  \ar@{^{(}->}[r]  \ar[rd]^{\varphi}&
   \PP({\mathcal E}^*) \ar[d]^{\mu} \\
& \PP^{n-1}  
}
}
o\`u $\varphi~: Z \to \PP^{n-1}$ est une contraction 
\'el\'ementaire, fibration en coniques,  
poss\'edant au moins une fibre singuli\`ere et o\`u 
${\mathcal E} = \varphi_*(-K_{Z})$.
Supposons que $D$ est un diviseur de $Z$ isomorphe \`a
$\PP^{n-1}$ de fibr\'e normal ${\mathcal O}_{\PP ^{n-1}}(d)$.
Alors $d$ est pair.
\end{lemm}
 
\noindent {\bf D\'emonstration.}
Soit $C$ une courbe de $Z$, composante d'une fibre singuli\`ere,
telle
que $-K_{\tilde Z} \cdot C =1$ (en particulier, $C$ est une droite
de $\PP({\mathcal E})^*$) et   
soit $a = D \cdot C$. Alors l'intersection de $D$
avec la fibre g\'en\'erique de $\varphi$ vaut $2a$
et $\varphi _{|D}: D \to \PP^{n-1}$ 
est finie de degr\'e (topologique) \'egal \`a $2a$.
Soit $r$ le degr\'e alg\'ebrique de $\varphi _{|D}: D \to \PP^{n-1}$
(autrement dit, $\varphi ^* {\mathcal O}_{\PP ^{n-1}}(1)_{|D}=
{\mathcal O}_{\PP ^{n-1}}(r)$)~: par le th\'eor\`eme de Bezout,
$2a = r^{n-1}$ donc $r$ et $a$ sont pairs.
Comme $\rho(Z /\PP^{n-1}) =1$, l'application de restriction
$\Pic (\PP ({\mathcal E}^*)) \to \Pic (Z)$
est un isomorphisme et il existe donc 
$b\in \ZZ$ tel que 
${\mathcal O}(D) = 
{\mathcal O}_{\PP ({\mathcal E}^*)}(a)_{|Z} 
\otimes
\varphi^* {\mathcal O}_{\PP ^{n-1}}(b)$. 
Comme ${\mathcal O}_{\PP ({\mathcal E}^*)}(1)_{|Z} = -K_Z$,
la formule d'adjunction donne
${\mathcal O}_{\PP ({\mathcal E}^*)}(1)_{|D} \simeq 
{\mathcal O}_{\PP ^{n-1}}(n+d)$,
d'o\`u ${\mathcal O}(D)_{|D} \simeq 
{\mathcal O}_{\PP ^{n-1}}(a(n+d)+br)$ donc
$d = a(n+d)+br$ et $d$ est pair.\finpreuve

\medskip
  
L'id\'ee de la d\'emonstration du th\'eor\`eme \ref{main} 
est simple~: la proposition \ref{Step1} produit
une vari\'et\'e de Fano $X'$. Dans le cas b), on peut
appliquer le lemme \ref{existcont} \`a $X'$ et $D'$
et \'etudier la contraction extr\'emale obtenue $\psi~: X' \to X''$. 
Le fait que $X'$ provienne de $\tilde X$ va nous permettre 
de d\'ecrire compl\`etement la situation. Ceci est
l'objet de la proposition suivante~:

\begin{prop}\label{Step2} 
Soient $\tilde X$ une vari\'et\'e de Fano 
de dimension $n \geq 3$, $\tilde D \simeq \PP ^{n-1}$
un diviseur de $\tilde X$ avec $N_{\tilde D / \tilde X}= 
{\mathcal O}_{\PP ^{n-1}}(-1)$. 
Soient $\varphi ~: \tilde X \to X'$
une contraction extr\'emale birationnelle 
d'ar\^ete $R = \QQ ^+ [\tilde F]$
avec $\tilde F \cdot \tilde D > 0$ 
et $\psi ~: X' \to X''$
une contraction extr\'emale 
d'ar\^ete $R' = \QQ ^+ [F']$
avec $F' \cdot \varphi (D) > 0$.
Alors $\psi ~: X' \to X''$ est
de l'un des deux types suivants~:
\begin{enumerate}
\item [a)] une application constante et $X' \simeq \PP ^n$~;
\item [b)] une fibration lisse en $\PP^1$ sur $\PP ^{n-1}$ donn\'ee par
$\psi~: X' \simeq \PP ({\mathcal O}_{\PP ^{n-1}} \oplus
{\mathcal O}_{\PP ^{n-1}}(d')) \to X'' \simeq \PP ^{n-1}$ 
avec $0 \leq d' \leq n-1$ et $D' = \varphi (D)$ est le diviseur
$\PP ({\mathcal O}_{\PP ^{n-1}})$.
\end{enumerate}
\end{prop}

\noindent {\bf D\'emonstration.}
Puisque $\varphi$ est birationnelle, nous sommes
dans la situation b) de la proposition \ref{Step1} 
dont nous gardons les notations~: 
$\varphi$ est une contraction lisse de centre $Y' \subset
\varphi (D) =D' \subset X'$ lisse 
de codimension deux.

\smallskip

\begin{enumerate}
\item [$\bullet$] Montrons tout d'abord que 
soit $\psi$ est constante (et donc $\rho (X') =1$),
soit les fibres non triviales de
$\psi$ sont de dimension un.
En effet, comme pour la  proposition \ref{Step1},
si $\psi$ poss\`ede une fibre de dimension au moins
deux, alors l'ar\^ete $R'$ contient une courbe rationnelle, not\'ee
$C$, contenue 
dans $D'$. Par cons\'equent, $d' =\deg (N_{D'/X'}) > 0$ et $C$ est une courbe
tr\`es libre, donc ses d\'eformations avec point fix\'e dans $D'$
couvrent un ouvert de $X'$, et par suite $\psi$ est constante
et $\rho(X') =1$. Dans ce cas, il d\'ecoule du lemme \ref{ro1}
ci-apr\`es que $X' \simeq \PP ^n$. 

\smallskip

\item [$\bullet$] Pla\c cons nous maintenant dans le cas o\`u 
les fibres non triviales de
$\psi$ sont de dimension un.
Comme pr\'ec\'edemment, par \cite{And85}, $X''$ est lisse et 
$\psi$ est de l'un des trois types suivants~: 
\begin{enumerate}
\item [(i)] une fibration lisse en $\PP ^1$,
\item [(ii)] une fibration en coniques,  
\item [(iii)] une contraction lisse de centre $Y'' \subset X''$ lisse 
de codimension deux.
\end{enumerate}

Nous allons montrer que seul le cas (i) est possible
et
l'observation suivante est cruciale~:
{\em soit $C' \subset X'$ une courbe non contenue 
dans $Y'$ (le centre de $\varphi$). Si
$-K_{X'}\cdot C' = 1$, alors $C'$ ne rencontre pas $Y'$
et si $-K_{X'}\cdot C' = 2$, alors $C'$ rencontre 
transversalement $Y'$ en exactement
un point.}
En effet, par l'absurde, si $\tilde C \subset \tilde X$ 
est la transform\'ee stricte de $C'$, alors  
$-K_{\tilde X} \cdot \tilde C = 
-K_{X'}\cdot C' - \tilde E \cdot \tilde C \leq 0$,
ce qui contredit $\tilde X$ Fano.

Excluons le cas (ii)~: comme l'union des fibres 
singuli\`eres est un diviseur de $X'$, son intersection avec 
$D'$ rencontre $Y'$ puisque $Y'$ est de codimension deux dans $X'$.  
Ainsi $Y'$ rencontre l'une des fibres scind\'ees $F' = F'_1 + F'_2$ de 
$\psi$, avec $-K_{X'}\cdot F'_i =1$, ce qui n'est pas possible d'apr\`es 
l'observation pr\'ec\'edente.

Excluons de m\^eme le cas (iii)~: les fibres $F'$ non triviales 
de $\psi$ v\'erifient $-K_{X'} \cdot F' =1$ et rencontrent toutes 
$D'$, donc certaines rencontrent $Y'$.

\smallskip

\item [$\bullet$] Etudions maintenant le cas (i)~: d'apr\`es \cite{Laz83},
$X''\simeq \PP ^{n-1}$ puisque $\psi _{|D'}~: D' \to X''$
est finie. 
Montrons que $\psi _{|D'}~: D' \to X''$ est un isomorphisme (de l\`a,
$D'$ est une section de $\psi$, et 
la conclusion b) s'ensuit !). Comme toutes les fibres $F'$
de $\psi$ v\'erifient $-K_{X'}\cdot F' = 2$, l'observation ci-dessus
montre que les fibres de $\psi$ rencontrant $Y'$ le font 
transversalement en exactement un point.
Il s'ensuit que $\psi _{|Y'}~: Y' \to Y'':= \psi (Y')$ est un isomorphisme.
Comme $Y'$ (resp. $Y''$) est une hypersurface
de l'espace projectif $D'$ (resp. $X''$), le th\'eor\`eme
de Lefschetz assure que $\psi ^*~: \Pic(X'') \to \Pic (D')$ est 
un isomorphisme et par suite que $\psi_{|D'}$ est un isomorphisme.  

\end{enumerate} 
\finpreuve
 
\medskip

\begin{lemm}\label{ro1} Soit $X'$ une vari\'et\'e projective 
lisse de dimension $n \geq 2$ contenant un diviseur $D' \simeq \PP ^{n-1}$.
Si $\rho (X') =1$, alors $X' \simeq \PP ^n$ et
$D'$ est un hyperplan.
\end{lemm}

\noindent {\bf D\'emonstration.} Puisque $\rho (X') =1$,
le diviseur $D'$ est ample et $X'$ est de Fano.
Comme les vari\'et\'es de Fano sont classifi\'ees en dimension
au plus trois o\`u le lemme se v\'erifie directement,
on peut supposer que $n \geq 4$. Le th\'eor\`eme de Lefschetz
assure alors que l'application de restriction
$\Pic (X') \to \Pic (D') =\ZZ \cdot {\mathcal O}_{\PP ^{n-1}}(1)$
est un isomorphisme. Soit donc ${\mathcal O}_{X'}(1)$
le g\'en\'erateur ample de $\Pic (X')$ dont la restriction 
\`a $D'$ est ${\mathcal O}_{\PP ^{n-1}}(1)$.
Soit $\nu$ l'indice de $-K_{X'}$. 
Par la formule d'adjonction, $-K_{X'|D'}= {\mathcal O}_{\PP ^{n-1}}(\nu)
= {\mathcal O}_{\PP ^{n-1}}(n+m)$ o\`u $m$ est le degr\'e
de $N_{D'/X'}$. Comme $m > 0$, il s'ensuit que $\nu \geq n+1$,
d'o\`u, d'apr\`es \cite{KO73}, $\nu = n+1$, $X' \simeq \PP ^n$
et $D'$ est un hyperplan.
\finpreuve
 
\section{D\'emonstration du th\'eor\`eme}

Soit $X$ une vari\'et\'e complexe de dimension $n \geq 3$
et supposons qu'il existe $a \in X$ tel que 
l'\'eclat\'ee de $X$ en $a$, not\'ee $\tilde X$, soit de Fano.
Notons $\tilde D$ le diviseur exceptionnel
de l'\'eclatement $\pi _a ~: \tilde X \to X$ et appliquons la 
proposition \ref{Step1}.
\begin{enumerate}
\item [$\bullet$] Dans le cas a),  
$\varphi~:\tilde X \simeq \PP ({\mathcal O}_{\PP ^{n-1}} \oplus
{\mathcal O}_{\PP ^{n-1}}(-1)) \to X' \simeq \PP ^{n-1}$ 
et $\tilde D$ est le diviseur
$\PP ({\mathcal O}_{\PP ^{n-1}})$. Il s'ensuit que $X$ est 
isomorphe \`a $\PP ^n$.
\item [$\bullet$] Dans le cas b), appliquons la proposition
\ref{Step2} et distinguons \`a nouveau deux cas~:
\begin{enumerate}
\item [$\bullet \bullet$] Dans le cas a), $X' \simeq \PP ^n$
et $\varphi~: \tilde X \to X' \simeq \PP ^n$ est l'\'eclatement de $\PP ^n$  
le long d'une quadrique ${\mathcal Q}_{n-2}$ contenue 
dans un hyperplan. Il s'ensuit que $X$ est la quadrique ${\mathcal Q}_{n}$
de dimension $n$. 
\item [$\bullet \bullet$] Dans le cas b), il y a un diagramme
commutatif~:
\medskip

\centerline{
\xymatrix{ & \tilde{X} \ar[ld]_{\pi_a}  \ar[rd]^{\varphi} \ar[d]_{\varphi '} 
& \\
X \ar[d]^{\mu} & \PP ({\mathcal O}
\oplus {\mathcal O}(1)) \ar[ld]_{\nu} \ar[rd]^{\psi '} &
X' \simeq \PP ({\mathcal O}
\oplus {\mathcal O}(d')) \ar[d]_{\psi} \\ 
\PP ^n & & X'' \simeq \PP^{n-1} \\  
}
} 
\medskip

\noindent o\`u $\mu$ est l'\'eclatement de $\PP ^n$ le 
long d'une sous-vari\'et\'e $Y$ de codimension deux dans $\PP^n$,
contenue dans un hyperplan $H$ avec $\mu (a) \notin H$
et $Y$ 
de degr\'e $d'+1$ dans $H$.

Explicitons en d\'etail le diagramme pr\'ec\'edent~:
$\varphi$ est l'\'eclatement de centre 
$Y' \subset D' = \PP({\mathcal O}_{\PP ^{n-1}})
\subset X' $ o\`u $Y'$ est un diviseur de $D'$ 
de degr\'e $d'+1$. Soit $E' = \psi^{-1} (\psi (Y'))$ l'union des
fibres de $\psi$ rencontrant $Y'$. Alors la transform\'ee stricte
$\tilde{E'}$ de $E'$ par $\varphi$ est exceptionnelle (voir par
exemple \cite{Mar82})
et correspond \`a l'\'eclatement $\varphi '$ de 
$\PP ({\mathcal O}_{\PP ^{n-1}}
\oplus {\mathcal O}_{\PP ^{n-1}}(1))$ 
de centre $Z \subset \PP({\mathcal O}_{\PP ^{n-1}})
\subset \PP ({\mathcal O}_{\PP ^{n-1}}
\oplus {\mathcal O}_{\PP ^{n-1}}(1))$ o\`u $Z$ est un diviseur 
de $\PP({\mathcal O}_{\PP ^{n-1}})$ de degr\'e $d'+1$.
Cette partie du diagramme est une transformation \'el\'ementaire
de Maruyama \cite{Mar82}.
L'application $\nu$ est la contraction sur un point
du diviseur ${\mathcal O}_{\PP ^{n-1}}(1) \subset
\PP ({\mathcal O}_{\PP ^{n-1}}
\oplus {\mathcal O}_{\PP ^{n-1}}(1))$
et $\mu$ est l'\'eclatement de $\PP ^n$
de centre $Y :=\nu(\varphi '(\tilde D \cap \tilde{E'}))$,
qui est un diviseur de 
degr\'e $d'+1$ de $H := \nu (\PP ({\mathcal O}_{\PP ^{n-1}}))
= \nu (\varphi ' (\tilde D))$. Enfin, $a = \mu ^{-1} ( 
\nu( \PP ({\mathcal O}_{\PP ^{n-1}}(1)))$.  
Autrement dit, $X \simeq V_{d'+1}$ avec les notations du th\'eor\`eme.
\end{enumerate} 
\end{enumerate} 
\finpreuve

-----------

\noindent L.B.~: {\em Institut Fourier, UFR de Math\'ematiques,
Universit\'e de Grenoble 1,
UMR 5582,
BP 74,
38402 Saint Martin d'H\`eres,
FRANCE. 

e-mail : bonavero@ujf-grenoble.fr
}

\noindent F.C.~: {\em Institut \'Elie Cartan (Math\'ematiques), 
Universit\'e H. Poincar\'e Nancy 1, UMR 7502, 
BP 239, 
54506 Vandoeuvre-l\`es-Nancy Cedex, 
FRANCE. 

e-mail : campana@iecn.u-nancy.fr
} 

\noindent J.A.W.~: {\em Instytut Matematyki UW
Banacha 2, PL-02-097 Warszawa, POLOGNE. 

e-mail : jarekw@mimuw.edu.pl
}
\end{document}